\documentstyle[leqno,hyperref]{book}

\newcounter{supersection}[section]
\newtheorem{th}[supersection]{Theorem}

\newtheorem{re}[supersection]{Remark}
\newtheorem{pr}[supersection]{Proposition}

\def\bibname{\textbf{REFERENCES}}
\def\thebibliography#1{\paragraph*{\uppercase{\bibname}}\list
{[\arabic{enumi}]}{\settowidth\labelwidth{[#1]}\leftmargin\labelwidth
\advance\leftmargin\labelsep\usecounter{enumi}}
\def\newblock{\hskip .11em plus .33em minus .07em}
\sloppy\clubpenalty4000\widowpenalty4000
\sfcode`\.=1000\relax}

\newfont{\newit}{cmfi10 scaled 1200}
\newfont{\newsit}{cmfi10 scaled 1000}
\newfont{\newsmit}{cmfi10 scaled 800}
\newfont{\newssmit}{cmfi10 scaled 600}

\def\F11{{}_{1}\mbox{\rm F}{}_{1}}

\begin{document}

\thispagestyle{plain}

\centerline{\Large \bf AN APPLICATION OF {\boldmath $\lambda$}-METHOD}

\medskip

\centerline{\Large \bf ON INEQUALITIES OF SHAFER-FINK'S TYPE}
\footnotetext{Research partially supported by the MNTRS, Serbia, Grant No. 144020.}

\vspace*{6.00 mm}

\centerline{\large \it Branko J. Male\v sevi\' c}

\vspace*{4.00 mm}

\begin{center}
\parbox{25.0cc}{\scriptsize \textbf{Abstract. In this article $\lambda$-method of Mitrinovi\' c-Vasi\' c
\cite{MitrinovicVasic70} is applied to improve the upper bound for the arc$\,$sin function of L.~Zhu \cite{Zhu05}.}}
\end{center}

\noindent
\section{\Large \bf \boldmath \hspace*{-7.0 mm}
1. Inequalities of \textbf{\textsc{Shafer}}-\textbf{\textsc{Fink}}'s type}


\noindent
{\rm D.$\,$S. Mitrinovi\' c} in \cite{MitrinovicVasic70} considered the lower bound
of the arc$\,$sin function, which belongs to {\rm R.$\,$E. Shafer}. Namely,
the following statement is true.

\begin{th}
For $0 \leq x \leq 1$ the following inequalities are true$:$
\begin{equation}
\displaystyle\frac{3x}{2 + \sqrt{1-x^2}}
\leq
\displaystyle\frac{6(\sqrt{1+x} - \sqrt{1-x})}{4 + \sqrt{1+x} + \sqrt{1-x}}
\leq
\mbox{\rm arc$\,$sin} \, x \, .
\end{equation}
\end{th}

\medskip
\noindent
{\rm A.$\,$M. Fink} proved the following statement in \cite{Fink95} .

\begin{th}
For $0 \leq x \leq 1$ the following inequalities are true$:$
\begin{equation}
\label{Ineq_Fink95}
\displaystyle\frac{3x}{2 + \sqrt{1-x^2}}
\leq
\mbox{\rm arc$\,$sin} \, x
\leq
\displaystyle\frac{\pi x}{2 + \sqrt{1-x^2}}.
\end{equation}
\end{th}

\medskip
\noindent
{\rm B.$\,$J. Male\v sevi\' c} proved the following statement in \cite{Malesevic97}.

\begin{th}
For $0 \leq x \leq 1$  the following inequalities are true$:$
\begin{equation}
\label{Ineq_Malesevic97}
\displaystyle\frac{3x}{2 + \sqrt{1-x^2}}
\leq
\mbox{\rm arc$\,$sin} \, x
\leq
\displaystyle\frac{\mbox{\small $\displaystyle\frac{\pi}{\pi-2}$} x}{
\mbox{\small $\displaystyle\frac{2}{\pi-2}$} + \sqrt{1-x^2}}
\leq
\displaystyle\frac{\pi x}{2 + \sqrt{1-x^2}}.
\end{equation}
\end{th}

\medskip
\noindent
The main result of the article \cite{Malesevic97} can be formulated with the next statement.

\begin{pr}
\label{Prop_01} In the family of the functions$:$
\begin{equation}
\label{Prop_Family_01}
f_{b}(x) = \displaystyle\frac{(b+1)x}{b + \sqrt{1-x^2}}
\quad (0 \leq x \leq 1),
\end{equation}
according to the parameter $b > 0$, the function $f_{2}(x)$ is the greatest lower bound of the
$\mbox{arc$\,$sin}\,x$ function and the function $f_{2/(\pi-2)}(x)$ is the least upper bound of the
$\mbox{arc$\,$sin}\,x$ function.
\end{pr}

\medskip
\noindent
{\rm L. Zhu} proved the following statement in \cite{Zhu05}.

\begin{th}
For $x \in [0,1]$  the following inequalities are true$:$
\begin{equation}
\label{Ineq_Zhu05}
\begin{array}{rcl}
\displaystyle\frac{3x}{2 + \sqrt{1-x^2}}
&\!\!\!\leq\!\!\!&
\displaystyle\frac{6(\sqrt{1+x} - \sqrt{1-x})}{4 + \sqrt{1+x} + \sqrt{1-x}}
\,\leq\,
\mbox{\rm arc$\,$sin} \, x                                                       \\[3.0 ex]
&\!\!\!\leq\!\!\!&
\displaystyle\frac{\mbox{\small $\pi(\sqrt{2}+\displaystyle\frac{1}{2})$}
(\sqrt{1+x} - \sqrt{1-x})}{4 + \sqrt{1+x} + \sqrt{1-x}}
\,\leq\,
\displaystyle\frac{\pi x}{2 + \sqrt{1-x^2}}.
\end{array}
\end{equation}
\end{th}

\medskip
\noindent
In this article we further improve the upper bound of the arc$\,$sin function. Namely, in the next
section we will give proof of the following theorem:

\begin{th}
\label{Th_second}
For $x \in [0,1]$ the following inequalities are true$:$
\begin{equation}
\label{Ineq_Malesevic06}
\begin{array}{rcl}
\displaystyle\frac{3x}{2 + \sqrt{1-x^2}}
&\!\!\!\leq\!\!\!&
\displaystyle\frac{6(\sqrt{1+x} - \sqrt{1-x})}{4 + \sqrt{1+x} + \sqrt{1-x}}
\,\leq\,
\mbox{\rm arc$\,$sin} \, x                                                       \\[3.0 ex]
&\!\!\!\leq\!\!\!&
\displaystyle\frac{\mbox{\small $\displaystyle\frac{\pi(2-\sqrt{2})}{\pi-2\sqrt{2}}$}( \sqrt{1+x} - \sqrt{1-x} )}{
\mbox{\small $\displaystyle\frac{\sqrt{2}(4-\pi)}{\pi-2\sqrt{2}}$} + \sqrt{1+x} + \sqrt{1-x}}
                                                                                 \\[5.0 ex]
&\!\!\!\leq\!\!\!&
\displaystyle\frac{\mbox{\small $\pi(\sqrt{2}+\displaystyle\frac{1}{2})$}
(\sqrt{1+x} - \sqrt{1-x})}{4 + \sqrt{1+x} + \sqrt{1-x}}
\,\leq\,
\displaystyle\frac{\pi x}{2 + \sqrt{1-x^2}}.
\end{array}
\end{equation}
\end{th}

\begin{re}
Using numerical method from {\rm \cite{Malesevic06c}} we have the following conclusions$:$

\medskip
\noindent
{\boldmath $1^{0}.$} For values $x \in (0,0.387 \, 266 \, 274 \ldots)$ the following inequality is true$:$
\begin{equation}
\mbox{\rm arc$\,$sin}\,x < \displaystyle\frac{\mbox{\small
$\displaystyle\frac{\pi}{\pi-2}$} x}{ \mbox{\small
$\displaystyle\frac{2}{\pi-2}$} + \sqrt{1-x^2}} <
\displaystyle\frac{\mbox{\small
$\pi(\sqrt{2}+\displaystyle\frac{1}{2})$} (\sqrt{1+x} -
\sqrt{1-x})}{4 + \sqrt{1+x} + \sqrt{1-x}},
\end{equation}
and for values $x \in (0.387 \, 266 \, 274 \ldots,1)$ the following inequality is true$:$
\begin{equation}
\mbox{\rm arc$\,$sin}\,x < \displaystyle\frac{\mbox{\small
$\pi(\sqrt{2}+\displaystyle\frac{1}{2})$} (\sqrt{1+x} -
\sqrt{1-x})}{4 + \sqrt{1+x} + \sqrt{1-x}} <
\displaystyle\frac{\mbox{\small $\displaystyle\frac{\pi}{\pi-2}$}
x}{ \mbox{\small $\displaystyle\frac{2}{\pi-2}$} + \sqrt{1-x^2}}.
\end{equation}
Numerically determined constant $c = 0.387 \, 266 \, 274 \ldots$ is the unique number where
the previous bounds have the same values over $(0,1)$.

\medskip
\noindent {\boldmath $2^{0}.$} For values $x \in (0, 1)$ the following inequality is true$:$
\begin{equation}
\mbox{\rm arc$\,$sin}\,x < \displaystyle\frac{\mbox{\small
$\displaystyle\frac{\pi(2-\sqrt{2})}{\pi-2\sqrt{2}}$}( \sqrt{1+x} -
\sqrt{1-x} )}{ \mbox{\small
$\displaystyle\frac{\sqrt{2}(4-\pi)}{\pi-2\sqrt{2}}$} + \sqrt{1+x} +
\sqrt{1-x}} < \displaystyle\frac{\mbox{\small
$\displaystyle\frac{\pi}{\pi-2}$} x}{ \mbox{\small
$\displaystyle\frac{2}{\pi-2}$} + \sqrt{1-x^2}}.
\end{equation}
\end{re}

\break

\noindent
\section{\Large \bf \boldmath \hspace*{-7.0 mm}
2. The main results} 

\noindent
In this article, using $\lambda$-method of Mitrinovi\' c-Vasi\' c we give an analogous statement
to Proposition \ref{Prop_01}. Let us notice that from inequality given by {\rm L. Zhu} \cite{Zhu05}:
\begin{equation}
\displaystyle\frac{6(\sqrt{1+x} - \sqrt{1-x})}{4 + \sqrt{1+x} +
\sqrt{1-x}} \,\leq\, \mbox{\rm arc$\,$sin} \, x \,\leq\,
\displaystyle\frac{\mbox{\small
$\pi(\sqrt{2}+\displaystyle\frac{1}{2})$} (\sqrt{1+x} -
\sqrt{1-x})}{4 + \sqrt{1+x} + \sqrt{1-x}},
\end{equation}
for $x \in [0,1]$, we can conclude that the function $\varphi(x) = \mbox{arc$\,$sin}\,x$ has a lower bound
and upper bound in the family of the functions:
\begin{equation}
\Phi_{\alpha,\beta}(x)
=
\displaystyle\frac{\alpha(\sqrt{1+x} - \sqrt{1-x})}{\beta + \sqrt{1+x} + \sqrt{1-x}}
\quad (0 \leq x \leq 1),
\end{equation}
for some values of parameters $\alpha, \beta > 0$. Next for $x = 0$ it is true that $\Phi_{\alpha,\beta}(0) =0$,
for $\alpha, \beta > 0$. On the other hand, for values $x \in (0, 1]$ it is true:
\begin{equation}
\label{Ineq_Comp_1}
\Phi_{\alpha_{1},\beta_{1}}(x) >
\Phi_{\alpha_{2},\beta_{2}}(x) \Longleftrightarrow \alpha_{1}
\beta_{2} - \alpha_{2} \beta_{1}
>
(\alpha_{2} - \alpha_{1}) (\sqrt{1\!+\!x} + \sqrt{1\!-\!x}),
\end{equation}
for $\alpha_{1,2}, \beta_{1,2} > 0$. Let us apply $\lambda$-method of Mitrinovi\' c-Vasi\' c on the considered
two-parameters family $\Phi_{\alpha,\beta}(x)$ in order to determine the bounds of the function $\varphi(x)$
under the following conditions:
\begin{equation}
\label{Cond_Ph}
\Phi_{\alpha,\beta}(0) = \varphi(0)
\quad \mbox{and} \quad
\displaystyle\frac{d}{dx} \Phi_{\alpha,\beta}(0) = \displaystyle\frac{d}{dx} \varphi(0).
\end{equation}
It follows that $\alpha = \beta + 2$. In that way we get one-parameter subfamily:
\begin{equation}
\label{Family_b}
f_{\beta}(x) = \Phi_{\beta+2,\beta}(x) =
\displaystyle\frac{(\beta+2)(\sqrt{1+x} - \sqrt{1-x})}{\beta +
\sqrt{1+x} + \sqrt{1-x}} \qquad (0 \leq x \leq 1),
\end{equation}
according to the parameter $\beta > 0$. For that family the condition (\ref{Cond_Ph}) is true:
\begin{equation}
\label{Cond_f1}
f_{\beta}(0) = \varphi(0)
\;\; \mbox{and} \;\;
\displaystyle\frac{d}{dx} f_{\beta}(0) = \displaystyle\frac{d}{dx} \varphi(0).
\end{equation}
Additionally, we have:
\begin{equation}
\label{Cond_f2}
\displaystyle\frac{d^{2}}{dx^{2}} f_{\beta}(0)
=
\displaystyle\frac{d^{2}}{dx^{2}} \varphi(0)
\;\; \mbox{and} \;\;
\displaystyle\frac{d^{3}}{dx^{3}} f_{\beta}(0)
=
\displaystyle\frac{d^{3}}{dx^{3}} \varphi(0) + \displaystyle\frac{4 \!-\! \beta}{4{\big (}2 \!+\! \beta{\big )}}
\end{equation}
and
\begin{equation}
\label{Cond_f3}
\;\;
\displaystyle\frac{d^{4}}{dx^{4}} f_{\beta}(0)
=
\displaystyle\frac{d^{4}}{dx^{4}} \varphi(0)
\;\; \mbox{and} \;\;
\displaystyle\frac{d^{5}}{dx^{5}} f_{\beta}(0)
=
\displaystyle\frac{d^{5}}{dx^{5}} \varphi(0) + \displaystyle\frac{3{\big (}128 \!+\! 18 \beta \!-\! 13\beta^2{\big )}}{
16 {\big (}2 \!+\! \beta{\big )}^2}.
\end{equation}
Let us notice that for the family of the functions $f_{\beta}(x)$, on the basis of (\ref{Ineq_Comp_1}),
for values $x \in (0,1]$ the following equivalence is true:
\begin{equation}
\label{Ineq_Comp_2} f_{\beta_{1}}(x) > f_{\beta_{2}}(x)
\Longleftrightarrow \beta_{1} < \beta_{2},
\end{equation}
for $\beta_{1,2} > 0$. Let us emphasize  that there is a better upper bound $f_{b_{1}}(x)$
than upper bound $\Phi_{\pi(\sqrt{2}+1/2), 4}(x)$ of the function $\varphi(x)$ over $(0,1]$.
It is true that the parameter $\beta = b_{1}$ fulfils:
\begin{equation}
f_{b_{1}}(1) = \varphi(1) = \displaystyle\frac{\pi}{2},
\end{equation}
hence:
\begin{equation}
b_{1} = \displaystyle\frac{\sqrt{2}(4 - \pi)}{\pi - 2\sqrt{2}} =
3.876 \, 452 \, 527 \, \ldots \, < 4 .
\end{equation}
Let us prove that the function $f_{b_{1}}(x)$ is the upper bound of the function $\varphi(x)$ over $[0,1]$.
Let us define the function:
\begin{equation}
h(x)
=
f_{b_{1}}(x) - \varphi(x)
\end{equation}
for $0 \leq x \leq 1$. For the function $h(x)$ we introduce two substitutions
$x = \cos t$ ${\big (}t \!\in\! [0, \mbox{\small $\displaystyle\frac{\pi}{2}$}]{\big )}$ and
$t = 4 \, \mbox{\rm arc}\,\mbox{\rm tg} \, u$ ${\big (}u \!\in\! [0, \mbox{\rm tg} \mbox{\small
$\displaystyle\frac{\pi}{8}$}]{\big )}$ respectively, and we get a new function:
\begin{equation}
\;\;
\mbox{\newsit{w}}(u)
\!=\!
h{\big (}\mbox{\rm cos}(4 \, \mbox{\rm arc}\,\mbox{\rm tg} \, u){\big )}
\!=\!
\displaystyle\frac{\sqrt{2}(b_{1}\!+\!2)(u^2\!+\!2u\!-\!1)}{
(\sqrt{2}\!-\!b_{1})u^2\!-\!2\sqrt{2}u\!-\!b_{1}\!-\!\sqrt{2}}
-
\displaystyle\frac{\pi}{2}
+
4 \, \mbox{\rm arc}\,\mbox{\rm tg} \, u
\end{equation}
for $0 \leq u \leq \mbox{\rm tg} \mbox{\small $\displaystyle\frac{\pi}{8}$} = \sqrt{2}\!-\!1$. Then:
\begin{equation}
\quad
\mbox{\small $\begin{array}{rcl}
\mbox{\normalsize $\displaystyle\frac{d}{d u} \mbox{\newsit{w}}(u)$}
\!&\!\!\!=\!\!\!&\!
{\Big (}
{\big (} \mbox{\normalsize $4 b_{1}^2$} \!+\! \mbox{\normalsize $2\sqrt{2}b_{1}^2$} \!-\! \mbox{\normalsize $8 b_{1}$}
\!-\! \mbox{\normalsize $4 \sqrt{2} b_{1}$} \!-\! \mbox{\normalsize $8$} {\big )} \mbox{\normalsize $u^4$} \!+\!
{\big (} \! \!-\!\mbox{\normalsize $4 \sqrt{2} b_{1}^2$} \!+\! \mbox{\normalsize $8 \sqrt{2} b_{1}$}
\!-\! \mbox{\normalsize $32$} {\big )} \mbox{\normalsize $u^3$} \!+\!                                         \\[1.5 ex]
\!&\!\!\! \!\!\!&\!
\;\;
{\big (} \mbox{\normalsize $8 b_{1}^2$} \!-\! \mbox{\normalsize $16 b_{1}$} \!-\! \mbox{\normalsize $16$} {\big )}
\mbox{\normalsize $u^2$} \!+\!
{\big (} \! \!-\!\mbox{\normalsize $4 \sqrt{2} b_{1}^2$} \!+\! \mbox{\normalsize $8 \sqrt{2} b_{1}$}
\!+\! \mbox{\normalsize $32$} {\big )} \mbox{\normalsize $u$}  \, +                                           \\[1.5 ex]
\!&\!\!\! \!\!\!&\!
\;\;
{\big (} \mbox{\normalsize $4 b_{1}^2$} \!-\! \mbox{\normalsize $2 \sqrt{2} b_{1}^2$} \!-\! \mbox{\normalsize $8 b_{1}$}
\!+\! \mbox{\normalsize $4 \sqrt{2} b_{1}$} \!-\! \mbox{\normalsize $8$} {\big )} {\Big )} {\Big /}            \\[1.5 ex]
\!&\!\!\! \!\!\!&\!
\;\;
{\Big (}
{\big (} \mbox{\normalsize $u^2$} \!+\! \mbox{\normalsize $1$} {\big )}
{\big (} \mbox{\normalsize $b_{1} u^2$} \!-\! \mbox{\normalsize $\sqrt{2} u^2$}
\!+\! \mbox{\normalsize $2 \sqrt{2} u$} \!+\! \mbox{\normalsize $b_{1}$} \!+\! \mbox{\normalsize $\sqrt{2}$} {\big )}^2 {\Big )}.
\end{array}$}
\end{equation}
All solutions of the equation $\displaystyle\frac{d}{d u} \mbox{\newsit{w}}(u) = 0$ are determined by terms:
\begin{equation}
\begin{array}{rcl}
u_{1,4} \!&\!\!=\!\!&\!
\displaystyle\frac{2 \sqrt{2} \mp \sqrt{-b_{1}^4 + 4 b_{1}^3 + 4 b_{1}^2 - 16 b_{1}}}{
b_{1}^2 - 2 b_{1} + 2 \sqrt{2} - 4 } \, , \\[1.5 ex]
u_{2,3} \!&\!\!=\!\!&\!
\sqrt{2} - 1 \,                         ;
\end{array}
\end{equation}
or by numerical values: $u_{1} \!=\! 0.0869 \ldots \,$, $u_{2,3} \!=\! 0.4142 \ldots \,$,
$u_{4} \!=\! 0.8400 \ldots \;$. The function $\mbox{\newsit{w}}(u)$ has local maximum at the point $u_{1}$ and
$\mbox{\newsit{w}}(0) \!=\! \mbox{\newsit{w}}(\sqrt{2} \!-\! 1) \!=\! 0$. Hence $\mbox{\newsit{w}}(u) \geq 0$
for $u \in [0, \sqrt{2} \!-\! 1]$. Therefore the function:
\begin{equation}
f_{b_{1}}(x)
=
\displaystyle\frac{\mbox{\small
$\displaystyle\frac{\pi(2-\sqrt{2})}{\pi-2\sqrt{2}}$}( \sqrt{1+x} -
\sqrt{1-x} )}{ \mbox{\small
$\displaystyle\frac{\sqrt{2}(4-\pi)}{\pi-2\sqrt{2}}$} + \sqrt{1+x} +
\sqrt{1-x}}
\end{equation}
is the upper bound of $\varphi(x)$ over $[0,1]$. Let us notice that, for values $x \in (0,1]$,
on the basis (\ref{Ineq_Comp_1}), the following inequalities are true:
\begin{equation}
\varphi(x)
<
f_{b_{1}}(x) = \Phi_{b_{1}+2,b_{1}}(x)
<
\Phi_{\pi(\sqrt{2}+1/2), 4}(x).
\end{equation}
Let us prove that the function $f_{b_{1}}(x)$ is the least upper bound of the function $\varphi(x)$ from
the family (\ref{Family_b}). The following implication is true:
\begin{equation}
\label{Impl_2} b_{1} < b \; \Longrightarrow \; f_{b}(1) <
f_{b_{1}}(1) = \varphi(1) = \displaystyle\frac{\pi}{2}.
\end{equation}
Hence for $b > b_{1}$ the function $f_{b}(x)$ is not the upper bound for the function $\varphi(x)$ over $[0,1]$.
According to the previous consideration we can conclude that the function $f_{b_{1}}(x)$ is the least upper bound
of the function $\varphi(x)$ over $[0,1]$.

\medskip
\noindent
The lower bound of the function $f_{4}(x)$ of the function $\varphi(x)$ over $[0,1]$, which belongs to
{\rm R.$\,$E. Shafer}, according to formulas (\ref{Cond_f1}) - (\ref{Cond_f3}), has at $x =0$ the root
of the fifth order. Let us prove that the function $f_{4}(x)$ is the greatest lower bound of the function
$\varphi(x)$ from the family (\ref{Family_b}). For fixed $b \in (b_{1},4)$ let us define the function:
\begin{equation}
g(x) = \left\{
\begin{array}{ccc}
\alpha                                                          & : & x = 0,       \\[2.0 ex]
\displaystyle\frac{f_{b}(x) - \varphi(x)}{x^3}                  & : & x \in (0,1];
\end{array}
\right.
\end{equation}
with the constant:
\begin{equation}
\alpha
=
\displaystyle\frac{\mbox{\small $\displaystyle\frac{d^3}{d x^3}$} f_{b}(0)
- \mbox{\small $\displaystyle\frac{d^3}{d x^3}$} \varphi(0)}{6}
=
\displaystyle\frac{4-b}{24{\big (}2+b{\big )}}
>
0.
\end{equation}
The function $g(x)$ is continuous over $[0,1]$ and the following is true:
\begin{equation}
g(0) > 0 \quad\mbox{and}\quad g(1) < 0.
\end{equation}
Therefore we can conclude that there is $c_{b} \in (0,1)$ such that $g(c_{b}) = 0$. Let us notice that
$g(0) > 0$ and $g(c_{b}) = 0$. Then, there is some point $\xi_{b} \in (0,c_{b})$ such that
$g(\xi_{b}) \!>\! 0$ {\big (}$g \!\in\! \mbox{\rm C}[0,c_{b}]${\big )}. This is sufficient
for conclusion that, for each $b \in (b_{1},4)$, the function $f_{b}(x)$ is not the lower bound of
the function $\varphi(x)$ over $[0,1]$. According to the previous consideration we can conclude
that the function $f_{4}(x)$ is the greatest lower bound of the function $\varphi(x)$ over $[0,1]$.

\medskip
\noindent
On the basis of the previous consideration the following statement is true.

\begin{pr}
\label{Prop_02} In the family of the functions$:$
\begin{equation}
\label{Prop_Family_02}
f_{b}(x) = \Phi_{b+2,b}(x) = \displaystyle\frac{(b+2)(\sqrt{1+x} -
\sqrt{1-x})}{b + \sqrt{1+x} + \sqrt{1-x}} \qquad (0 \leq x \leq 1),
\end{equation}
according to the parameter $b > 0$, the function $f_{4}(x)$ is the greatest lower bound of the
$\mbox{arc$\,$sin}\,x$ function and the function $f_{\sqrt{2}(4 - \pi)/(\pi - 2\sqrt{2})}(x)$ is
the least upper bound of the $\mbox{arc$\,$sin}\,x$ function.
\end{pr}

\begin{re}
Let us emphasize that Theorem {\rm \ref{Th_second}} has been recently considered
in \mbox{\rm \cite{Zhu07}}~and~\mbox{\rm \cite{Zhu07b}}. In the article
{\rm \cite{Zhu07b}} a simple proof of Theorem {\rm \ref{Th_second}}
based on "L'Hospital ru\-le for mono\-tonicity" is obtained.
\end{re}

\bigskip
{\small
\noindent University of Belgrade,
          \hfill (Received     09/30/2006)                                        \break
\noindent Faculty of Electrical Engineering,
          \hfill $($Revised 05/08/2007$)\,$                                       \break
\noindent P.O.Box 35-54, $11120$ Belgrade, Serbia                           \hfill\break
\noindent {\footnotesize \bf malesevic@etf.bg.ac.yu}, {\footnotesize \bf malesh@eunet.yu}
\hfill}

\end{document}